\documentclass[12pt]{amsart}

\usepackage{url}

\newtheorem{theorem}{Theorem}

\begin{document}

\title{Vita: Friedrich Wilhelm Wiener}

\author{Harold P. Boas}

\address{Department of Mathematics, Texas A\&M University,
  College Station, TX 77843--3368}

\email{boas@math.tamu.edu}

\thanks{The authors' research was partially supported by grants
  from the National Science Foundation.}

\author{Dmitry Khavinson}

\address{Department of Mathematical Sciences, University of
  Arkansas, Fayetteville, AR 72701} 

\email{dmitry@comp.uark.edu}

\begin{abstract}
  We report on the life and work of F.~Wiener, whom we confused
  with N.~Wiener in a previous article.
\end{abstract}


\maketitle

\section{Introduction}
In a recent note \cite{BoasKhavinson}, we proved a
multi-dimensional analogue of the following classical theorem of
Harald Bohr~\cite{Bohr}. (For subsequent developments in the
multi-dimensional theory, see \cite{Aizenberg-Bohr,
  Aizenberg-Aytuna-Djakov-generalization,
  Aizenberg-Aytuna-Djakov-abstract}.)

\begin{theorem}[Bohr]
  Suppose that a power series $\sum_{k=0}^\infty c_k z^k$
  converges for $z$ in the unit disk, and $|\sum_{k=0}^\infty c_k
  z^k|< 1$ when $|z|< 1$.  Then $\sum_{k=0}^\infty |c_k z^k|< 1$
  when $|z|< 1/3$. Moreover, the radius $1/3$ is the best
  possible.
\end{theorem}

In one part of the proof, we adapted to higher dimensions an
elegant argument that Bohr attributed to Wiener.  Since Bohr
mentioned this name in the same sentence with the names of Riesz
and Schur, we assumed it to be the famous Norbert Wiener, and we
added the initial~``N'' in our attribution. Our assumption was
false. Lawrence Zalcman brought to our attention that Edmund
Landau mentioned the name of one F.~Wiener in connection with
Bohr's theorem \cite[\S4]{LandauE}.

\section{Wiener's life}
Having never heard of a mathematician F.~Wiener, we investigated.
We report here on what information we have discovered about the
life and work of F.~Wiener, hoping that his name may be preserved
in mathematical history for another generation.

According to the curriculum vitae accompanying his dissertation,
Friedrich Wilhelm Wiener was born in 1884 in Meseritz, then part
of the Prussian province of Posen and now part of Poland.  After
completing high school (gymnasium), he pursued studies in
G\"ottingen. After a year of compulsory military service in
1904--1905, he resumed studies in Berlin. He returned to
G\"ottingen in 1909, the same year that Landau was called there
as Minkowski's successor. Wiener attended lectures of such famous
mathematicians as Frobenius, Hilbert, Landau, Schottky, Schur,
and Schwarz. He completed his doctoral dissertation
\cite{WienerF1911} under the supervision of Landau in 1911.

Wiener published one journal article \cite{WienerF1910} in 1910,
which is cited in standard books \cite{HardyLittlewoodPolya,
  PolyaSzego}.  After a promising beginning, he seems to have
published nothing further, not even his dissertation.  There is
no evidence that Wiener was ever a member of the Deutsche
Mathematiker-Vereinigung (DMV); no obituary notice for Wiener
appeared in the DMV \emph{Jahresbericht} \cite{DMVyearbook}.
Although we do not know the circumstances of Wiener's death, this
must have occurred no later than 1921, as the index published
that year to volumes 51--80 of \emph{Mathematische Annalen} lists
Wiener as deceased.  We conjecture that Wiener may have been a
casualty of the war.

\section{Wiener's work}
The focus of Wiener's mathematical work was to discover simple
proofs of known theorems. Both of his papers have the word
``elementary'' in the title. 

\subsection{Hilbert's inequality}
Wiener's 1910 paper concerns Hilbert's double series theorem
stating the boundedness in~\(\ell_2\) of the quadratic form \(
\sum_{m=1}^\infty \sum_{n=1}^\infty x_m x_n/(m+n)\).

\begin{theorem}[Hilbert]
  \begin{equation*}
    \biggl| \sum_{m=1}^\infty \sum_{n=1}^\infty
    \frac{x_mx_n}{m+n} \biggr| \le C \sum_{n=1}^\infty |x_n|^2. 
  \end{equation*}
Moreover, the inequality holds with \(C=\pi\), and no smaller
value of the constant~\(C\) will do.
\end{theorem}

Hilbert's proof was first published in the dissertation
\cite{Weyl} of his student Hermann Weyl in 1908.  The theorem
attracted a great deal of attention, and numerous proofs and
generalizations were published subsequently.  The classical book
by Hardy, Littlewood, and P\'olya \cite{HardyLittlewoodPolya}
devotes a whole chapter to this inequality. At the time of
Wiener's work, it was not known that the sharp value of the
constant~\(C\) is~\(\pi\): Schur proved this the following
year~\cite{Schur}.

What Wiener meant by an ``elementary'' proof of Hilbert's
inequality was a proof that used no integration and no function
theory. His proof consists of the following elementary steps.
\begin{enumerate}
\item Reduce to the case that \(\{x_n\}_{n=1}^\infty\) is a
  decreasing sequence of positive real numbers.
\item Group the terms in the inner sum into blocks whose terms
  have indices running between consecutive squares.
\item Apply the Cauchy-Schwarz inequality to both the inner sum
  and the outer sum.
\item Interchange the order of summation.
\item Invoke Cauchy's condensation test for convergence of
  series. 
\end{enumerate}

\subsection{Wiener's dissertation}
In his dissertation, Wiener addresses two questions in the theory
of entire functions of one complex variable.

The first part of the dissertation concerns the minimum modulus
of an entire function~\(f\). Let \(m(r)= \min\{|f(re^{i\theta})|:
0\le\theta\le2\pi\}\). Since \(m(r)\) is zero when \(f\)~has a
zero of modulus~\(r\), the natural question to ask about a lower
bound for \(m(r)\) is whether \(m(r)\) is frequently large: is
there some reasonable comparison function \(c(r)\) such that \(
\limsup_{r\to\infty} m(r)/c(r)>0\)?

If \(f\)~is an entire function of finite order at most~\(\rho\),
meaning that \( \lim_{|z|\to\infty}
|f(z)|e^{-|z|^{\rho+\epsilon}}=0\) for every
positive~\(\epsilon\), then Hadamard's factorization theorem
implies that \( \limsup_{r\to\infty}
m(r)e^{r^{\rho+\epsilon}}=\infty\) for every
positive~\(\epsilon\).  In other words, \(m(r)\) cannot tend to
zero too fast.  This weak estimate cannot be improved in general.
For example, the exponential function~\(e^z\) has order~\(1\) and
\(m(r)=e^{-r}\).  On the other hand, if \(f\)~is a non-constant
polynomial, then \(m(r)\) tends to infinity like a power
of~\(r\).  The question arises of whether an entire function of
sufficiently small order is enough like a polynomial that its
minimum modulus must be unbounded.

In 1905, A.~Wiman confirmed \cite{Wiman} that the minimum modulus
of every non-constant entire function of order~\(\rho\) strictly
less than \(1/2\) is indeed unbounded. Moreover, \(
\limsup_{r\to\infty} m(r) e^{-r^{\rho-\epsilon}}=\infty\) when
\(0<\epsilon<\rho<1/2\). The cutoff at \(1/2\) is sharp, for the
convergent infinite product \( \prod_{n=1}^\infty (
1-\frac{z}{n^2} )\), which equals \(
(\sin\pi\sqrt{z}\,)/(\pi\sqrt{z}\,)\), has order \(1/2\) and
\(m(r)\le r^{-1/2}/\pi\). (See \cite[Chapter~3]{Boas-entire} for
more about the minimum modulus of entire functions of small
order.)

Wiener's dissertation gives a new proof of Wiman's theorem. The
proof is elementary in the sense that it uses only arguments
about series and products of real numbers; it avoids using
theorems from function theory.

Wiener's proof even supplements Wiman's theorem by giving some
information in the endpoint cases \(\rho=0\) and \(\rho=1/2\).
Namely, Wiener shows that if \(f(z) = \prod_{n=1}^\infty (
1-\frac{z}{a_n}) \), where \(\{a_n\}_{n=1}^\infty\) is a sequence
of non-zero complex numbers of increasing modulus, and if
\(\lim_{n\to\infty} n^2/|a_n|=0\), then \( \limsup_{r\to\infty}
m(r)r^{-k}=\infty\) for every positive~\(k\). This result applies
to all transcendental entire functions of order~\(0\) (for
example, to \(\prod_{n=1}^\infty (1-\frac{z}{n^n})\)) and to some
entire functions of order~\(1/2\) (for example, to
\(\prod_{n=2}^\infty (1-\frac{z}{n^2\log n})\)).

The second part of Wiener's dissertation is motivated by a
theorem of Landau \cite{Landau} that generalizes Picard's little
theorem.

\begin{theorem}[Landau]
  There is a positive function \(R\) such that every polynomial
  of the form \(a_0+z+a_2 z^2+\dots+a_nz^n\) assumes at least one
  of the values \(0\) and~\(1\) in the disk \(\{z: |z|\le
  R(a_0)\}\). The function~\(R\) is independent of the
  degree~\(n\) and the higher coefficients \(a_2\), \dots,
  \(a_n\).
\end{theorem}

One might hope that a theorem about polynomials would have an
elementary proof, which would then yield an elementary proof of
Picard's theorem. Wiener was able to find an elementary proof
(using Rouch\'e's theorem, but nothing else from function theory)
of Landau's theorem under an additional hypothesis about the
location of the zeroes of the polynomial. Namely, he assumed that
the zeroes are located within the two equal acute angles
determined by two lines intersecting at the origin. If the radian
measure of the acute angle is \(\frac{1}{2}\pi - \beta\), then
one can take \(R(a_0) = 28|a_0 \log a_0|/\sin\beta\). (The cases
\(a_0=0\) and \(a_0=1\) are of no concern, because then the
polynomial takes the value \(0\) or~\(1\) at the origin.)

\section{Acknowledgments}
For assistance in this project of identifying and tracing
F.~Wiener, we thank Samuel~J. Patterson
(Georg-August-Universit\"at G\"ottingen), Constance Reid,
Heinrich Wefelscheid (Gerhard-Mercator-Universit\"at
Gesamthochschule Duisburg), and Lawrence~A. Zalcman (Bar Ilan
University).  We are especially indebted to Professor Wefelscheid
for locating and sending to us a copy of Wiener's dissertation.
We thank Heidemarie W\"ormann Boas for help with German
translation.

\bibliographystyle{amsplain} 
\bibliography{wiener}

\end{document}